\documentclass[12pt]{amsart}

\usepackage[a-2b,mathxmp]{pdfx}[2018/12/22]%for NSF report

\usepackage{newpxtext,newpxmath}

\usepackage[bbgreekl]{mathbbol}
\usepackage[shortlabels]{enumitem}
\usepackage{graphicx}

\usepackage{hyperref}
\hypersetup{
	pdfstartview={XYZ null null 1.00}, 
	pdfpagemode=UseNone, 
	colorlinks,
	breaklinks, 
	linkcolor=blue,
	urlcolor=blue, 
	anchorcolor=blue,
	citecolor=blue
}

\usepackage{comment}

\usepackage{ifthen}
\usepackage{mathrsfs}
\usepackage{pict2e}
\usepackage{xargs}
\usepackage{xspace}

\usepackage[capitalise]{cleveref}
\usepackage{geometry}
\DeclareSymbolFontAlphabet{\amsmathbb}{AMSb}

\newcommand{\definedterm}[1]{\textbf{#1}}

\let\e\varepsilon
\let\0\emptyset

\let\imp\Rightarrow

\newcommand{\absolutevalue}[1]{|#1|}
\newcommand{\Bairetree}[1][]{
\ifthenelse{\equal{#1}{}}{\functions{<\N}{\N}}{\functions{#1}{\N}}
}
\newcommand{\bbone}{\mathbb{1}}

\newcommand{\calB}{\mathcal{B}}

\newcommand{\calR}{\mathcal{R}}
\newcommand{\Cantorspace}[1][]{
\ifthenelse{\equal{#1}{}}{\functions{\N}{2}}{\functions{#1}{2}}
}
\newcommand{\Cantortree}[1][]{
\ifthenelse{\equal{#1}{}}{\functions{<\N}{2}}{\functions{#1}{2}}
}

\newcommand{\characteristicfunction}[1]{\bbone_{#1}}

\newcommand{\composition}{\circ}
\newcommandx{\concatenation}[2][1 = undefined, 2 = undefined]{
\ifthenelse{\equal{#1}{undefined}}{{}\smallfrown}{
\ifthenelse
  {\equal{#2}{undefined}}
  {\smallfrown_{#1}}
  {\bigoplus_{#1} #2}
}
}
\newcommandx{\disjointunion}[2][1 = undefined, 2 = undefined]{
\ifthenelse{\equal{#1}{undefined}}{\sqcup}{
\ifthenelse{\equal{#2}{undefined}}{\coprod #1}{\coprod_{#1} #2}
}
}

\newcommand{\forwardorbit}[2]{[#1]_{\overrightarrow{#2}}}
\newcommand{\backwardorbit}[2]{[#1]_{\overleftarrow{#2}}}
\newcommand{\from}{\colon}
\newcommandx{\functions}[3][3 =]{
\ifthenelse{\equal{#3}{}}{#2^{#1}}{#2^{#1}_{#3}}
}
\newcommand{\goesto}{\rightarrow}
\newcommand{\graph}[1]{\mathrm{graph}(#1)}
\newcommand{\heightcorrection}[1]{\raisebox{0pt}[0pt][0pt]{#1}}
\newcommand{\hittingtime}[2]{n_{#2}^{#1}}
\newcommand{\image}[3][]{
\ifthenelse{\equal{#1}{}}{#2(#3)}{#2^{#1}(#3)}
}

\newcommandx{\integral}[4][1=,4=]{
\ifthenelse{\equal{#1}{}}{\int #2}{\int_{#1} #2}
\ifthenelse{\equal{#2}{}}{d#3}{\ d#3}
\ifthenelse{\equal{#4}{}}{}{(#4)}
}

\newcommand{\Exp}[1]{\E(#1 \mid \calB_T)}

\newcommand{\dm}[1][\mu]{\, d #1}

\newcommandx{\intersection}[2][1 =, 2 =]{
\ifthenelse{\equal{#1}{}}{\cap}{
\ifthenelse{\equal{#2}{}}{\bigcap #1}{{\bigcap_{#1} #2}}
}
}

\newcommand{\lexeq}{\le_{\mathrm{lex}}}
\newcommand{\limit}[2]{\ratio{#1}{#2}{\infty}}
\newcommand{\limitpoints}[2]{\calR^\coc(#1, #2)}
\newcommand{\limitinferior}[2]{\underline{R}^\coc(#1, #2)}
\newcommand{\limitsuperior}[2]{\overline{R}^\coc(#1, #2)}
\newcommand{\longalign}{\hspace*{10pt} & \hspace*{-10pt}}

\newcommand{\N}{\amsmathbb{N}}

\newcommand{\periodicpart}[3][]{\text{Per}_{#1}^{#2}(#3)}
\newcommand{\positiveintegers}{\N^+}
\newcommand{\positivereals}{\R^+}

\newcommand{\preimage}[3][]{
\ifthenelse{\equal{#1}{}}{#2^{-1}(#3)}{#2^{-#1}(#3)}
}
\newcommandx{\product}[2][1 =, 2 =]{
\ifthenelse{\equal{#1}{}}{\times}{
\ifthenelse{\equal{#2}{}}{\prod #1}{{\prod_{#1} #2}}
}
}

\newcommand{\R}{\amsmathbb{R}}
\newcommand{\ratio}[3]{R^\coc_{#3}(#1, #2)}
\renewcommand{\restriction}[2]{#1 |_{#2}}
\newcommandx{\sequence}[2][2 = undefined]{
\ifthenelse{\equal{#2}{undefined}}{(#1)}{
(#1)_{#2}
}
}
\newcommandx{\set}[2][2 = undefined]{
\ifthenelse{\equal{#2}{undefined}}{\{ #1 \}}{
\{ #1 \suchthat #2 \}
}
}
\newcommand{\setcomplement}[1]{\twiddle #1}
\newcommand{\suchthat}{\mid}
\newcommandx{\summation}[2][1 =, 2 =]{
\ifthenelse{\equal{#1}{}}{+}{
\ifthenelse{\equal{#2}{}}{\sum #1}{{\sum_{#1} #2}}
}
}

\newcommand{\twiddle}
{\raisebox{1.5pt}{\scalebox{.75}{$\mathord{\sim}$}}}
\newcommandx{\union}[2][1 =, 2 =]{
\ifthenelse{\equal{#1}{}}{\cup}{
\ifthenelse{\equal{#2}{}}{\bigcup #1}{{\bigcup_{#1} #2}}
}
}
\newcommand{\uppersum}[2]{S^\coc_{#2} #1}

\newcommand{\E}{\amsmathbb{E}}

\let\~\overline

\newcommand*{\defeq}{\mathrel{\vcenter{\baselineskip0.5ex \lineskiplimit0pt \hbox{\scriptsize.}\hbox{\scriptsize.}}}=}

\newcommand{\coc}{w}

\newcommand{\Xas}{X_\mathrm{as}}
\newcommand{\Xep}{X_\mathrm{ep}}

\newcommand{\Dowker}{Dow\-ker\xspace}

\newenvironment{lemmaproof}{

\begin{proof}
}{\end{proof}}

\newenvironment{propositionproof}{

\begin{proof}
}{\end{proof}}

\newenvironment{theoremproof}{

\begin{proof}
}{\end{proof}}

\numberwithin{equation}{section}

\crefname{introtheorem}{Theorem}{Theorems}

\crefname{lemma}{Lemma}{Lemmas}
\newtheorem{lemma}[equation]{Lemma}

\crefname{proposition}{Proposition}{Propositions}
\newtheorem{proposition}[equation]{Proposition}

\crefname{theorem}{Theorem}{Theorems}
\newtheorem{theorem}[equation]{Theorem}

\crefname{corollary}{Corollary}{Corollaries}

\theoremstyle{definition}

\crefname{introremark}{Remark}{Remarks}

\crefname{remark}{Remark}{Remarks}

\newtheorem*{claim}{Claim}

\begin{document}

%% Front matter

\author[B.D. Miller]{Benjamin D. Miller}
\address{
Benjamin D. Miller \\
1008 Balsawood Drive \\
Durham, NC 27705 \\
USA
}
\email{glimmeffros@gmail.com}
\urladdr{
\url{http://sites.google.com/view/b-miller}
}

\author[A. Tserunyan]{Anush Tserunyan}
\address{
Anush Tserunyan \\
Mathematics and Statistics Department \\
McGill University \\
Montr\'{e}al, QC \\
Canada
}
\email{anush.tserunyan@mcgill.ca}
\urladdr{
\url{https://www.math.mcgill.ca/atserunyan}
}



\title{A ratio ergodic theorem via tiling and uniformly syndetic markers}

\begin{abstract}
We prove a purely Borel/measureless version of Dowker's ratio ergodic theorem, from which we derive a strengthening of Dowker's original theorem with a precise identification of the limit of local ergodic ratios.
This is done by implementing the pointwise tiling idea of \cite{Tserunyan} in the more complex setting of continuum-to-one Borel transformations.
Along the way, we establish a vanishing markers lemma for these transformations, which generalizes its well-known counterpart for invertible transformations.
\end{abstract}

\maketitle

\section{Introduction}

The classical pointwise ergodic theorem, which dates back to Birkhoff \cite{Birkhoff} and Khintchine \cite{Khintchine:Birkhoff_thm}, states that for every measure-preserving (not necessarily invertible) transformation $T$ of a probability space $(X, \mu)$ and every $f \in L^1(X, \mu)$, the averages $\frac{S_n f}{n}$, where $S_n f \defeq \sum_{k<n} f \circ T^k$, converge a.e.\ to the conditional expectation $\Exp{f}$ of $f$ with respect the $\sigma$-algebra $\calB_T$ of $T$-invariant Borel sets. 
Attempting to generalize this theorem to $\sigma$-finite measures, we see that if the measure is infinite and $T$ is ergodic then the averages $\frac{S_n f}{n}$ converge to $0$ a.e., so we have to at least give up the statement about conditional expectation.
Despite this, Hopf \cite[Sec.14, pp.49-50]{Hopf:book1937} proved that for a \textit{conservative invertible measure-preserving} transformation $T$ of a $\sigma$-finite measure space $(X,\mu)$ and functions $f \in L^1(X, \mu)$ and $g > 0$, the ratios $\frac{S_n f}{S_n g}$ converge a.e.\ to a finite limit.
(A particular case of this was proved a year earlier by Stepanoff \cite{Stepanoff:Hopf_thm}.)

Hopf's theorem was later generalized by Hurewicz \cite{Hurewicz} and Halmos \cite{Halmos:ergodic} to all \textit{invertible nonsingular} transformations, replacing the sums $S_n f$ and $S_n g$ with weighted sums $\uppersum{f}{n}$ and $\uppersum{g}{n}$, where the weight function $w$ is the Radon--Nikodym derivative $d(T^{-1}_* \mu)/d \mu$.
Finally, Dowker \cite[Theorem~II]{Dowker} removed the invertibility hypothesis, proving for \textit{all nonsingular} transformations that the limit of $\frac{\uppersum{f}{n}}{\uppersum{g}{n}}$ exists and is finite a.e.%
\footnote{In Dowker's theorem, the weight function $w$ is $\frac{d\mu}{d(T_* \mu)} \circ T$, which is equal to $\frac{d(T^{-1}_* \mu)}{d\mu}$ for invertible $T$.}
For integrable $g$, Dowker's theorem was ultimately generalized%
\footnote{Indeed, taking $P f \defeq (T f) \cdot w$ we obtain Dowker's theorem from the Chacon--Ornstein theorem.}
by the well-known Chacon--Ornstein theorem \cite{Chacon-Ornstein}, which states the a.e.\ existence and finiteness of the limit of $\frac{\sum_{i<n}P^i f}{\sum_{i<n} P^i g}$ for every \textit{positive linear operator} $P$ on $L^1$ of norm $\le 1$.

Our main result is a purely Borel/measureless version (\cref{intro:measureless}) of Dowker's theorem.
From this we derive a strengthening (\cref{intro:Dowker_with_limit}) of Dowker's theorem with the precise identification of the limit of $\frac{\uppersum{f}{n}}{\uppersum{g}{n}}$.
For conservative $T$, our identification of the limit coincides with that of Chacon \cite{Chacon:identification_of_limit} and Neveu \cite[Proposition~V.6.4]{Neveu:book}.

The proof of our measureless ergodic theorem and the derivation of the strengthening of Dowker's theorem from it uses a significant generalization of the pointwise tiling idea of the second author, given in \cite{Tserunyan} in a much simpler setting.
In particular, to prove the measureless version, we establish a vanishing markers lemma for (potentially) continuum-to-one Borel transformations (\cref{aperiodic:markers}), which generalizes and strengthens its well-known counterpart for invertible transformations \cite[\S3,~Lemma~1]{SlamanSteel}.

\subsection*{Statements of results}

We now establish necessary notation and terminology to state our strengthening of Dowker's theorem and the underlying purely Borel version.

\medskip

Let $T \from X \to X$ be a \definedterm{positively nonsingular} transformation of a $\sigma$-finite measure space $(X,\mu)$, i.e.\ $\mu \ll T_* \mu$. 
Define $\coc \from \N \times X \to \positivereals$ by $(x,k) \mapsto \frac{d\mu}{d(T^k_* \mu)}(T^k x)$. 
We often write $\coc_k(x)$ or $\coc^x(k)$ to mean $\coc(k,x)$.
By the chain rule,
\[
\coc_k(x) 
=
\tfrac{d\mu}{dT_* \mu}(T^k x)
\tfrac{d\mu}{dT_* \mu}(T^{k-1} x) 
\cdots
\tfrac{d\mu}{dT_* \mu}(T x) 
,
\]
so $\coc$ is a \definedterm{Borel $T$-cocycle} to the multiplicative group $\positivereals$ of positive reals, that is, a Borel map $\N \times X \to \positivereals$ satisfying the \definedterm{cocycle identity}
\begin{equation}\label{eq:coc_identity}
\coc_{k + \ell}(x) = \coc_k(T^\ell x) \cdot \coc_\ell(x)
\end{equation}
for all $x \in X$ and $k, \ell \in \N$.
It also follows from the definitions of the Radon--Nikodym derivative $\frac{d\mu}{d(T^k_* \mu)}$ and the pushforward measure $(T^k_* \mu)$ that $\mu$ is \definedterm{$T$-$\coc$-invariant}, i.e.
\begin{equation}\label{eq:rho-invariance}
\int f \dm 
=
\int f (T^k x) \coc_k(x) \dm(x),
\end{equation}
for each non-negative or $\mu$-integrable measurable function $f \from X \to \R$ and $k \in \N$. 
We refer to this function $\coc$ as the \definedterm{Radon--Nikodym cocycle} associated with $(X,\mu,T)$. 
Intuitively, $\coc_k(x)$ is equal to the ``weight'' of $T^k x$ relative to the ``weight'' of the set $T^{-k}(T^k x)$.

For $f \from X \to \R$, $n \in \N$, and $x \in X$, put
\[
\uppersum{f}{n}(x) \defeq \sum_{k<n} f(T^k x) \coc_k(x).
\]
Furthermore, for $g \from X \to \positivereals$, put 
\[
\ratio{f}{g}{n} \defeq \frac{\uppersum{f}{n}}{\uppersum{g}{n}} 
%\text{ and } A^\coc_n f(x) \defeq \ratio{f}{1}{n},
\]
and call it the ($\coc$-weighted) \definedterm{average-ratio} of $f$ and $g$ over the set $\set{x,Tx,\dots, T^{n-1}x}$.
%so $A^\coc_n f(x)$ is the \definedterm{$\coc$-weighted average} of $f$ over the set $\set{x,Tx,\dots, T^{n-1}x}$.

\begin{theorem}[Dowker \cite{Dowker}]\label{Dowker}
Let $T$ be a positively nonsingular transformation of a $\sigma$-finite measure space $(X,\mu)$ and let $\coc$ be the associated Radon--Nikodym cocycle. 
For each $f \in L^1(X,\mu)$ and measurable $g \from X \to \positivereals$ with $\lim_n \uppersum{g}{n} = \infty$ a.e., the limit
\[
\limit{f}{g} \defeq \lim_n \ratio{f}{g}{n}
\]
exists a.e.\ and is finite.
\end{theorem}

Recall that a set $Y \subseteq X$ is called \definedterm{$T$-invariant} if $T^{-1} Y = Y$.
We denote by $\Exp{h}$ the conditional expectation of a function $h \in L^1(X,\mu)$ with respect to the $\sigma$-algebra $\calB_T$ of all $T$-invariant Borel subsets of $X$.
Our strengthening then says that the limit is as expected:

\begin{theorem}[Ratio ergodic with limit identified]\label{intro:Dowker_with_limit}
    Under the hypothesis of \cref{Dowker}, we moreover have that the limit $\limit{f}{g}$ is $T$-invariant and:
    \begin{enumerate}[(i)]
    \item If $g$ is $\mu$-integrable, then $\limit{f}{g} = \frac{\Exp{f}}{\Exp{g}}$ a.e.

    \item If $\int_Y g \dm = \infty$ for each non-null $Y \in \calB_T$, then $\limit{f}{g} = 0$ a.e.
\end{enumerate}
\end{theorem}

\cref{intro:Dowker_with_limit} follows from a purely Borel counterpart. 
To state it we need to describe the Borel counterparts/witnesses of the measure-theoretic notions involved in \cref{Dowker}.

Firstly, we replace the measure space by a \definedterm
{Borel space}, i.e.\ a set $X$ equipped with a distinguished $\sigma$-algebra $\calB$ of subsets of $X$, which we call Borel%
\footnote{Note that $\calB$ may contain all $\mu$-measurable sets for some measure $\mu$ defined on a topological space $X$, so it may be a larger $\sigma$-algebra than that of topologically Borel subsets of $X$.}%
.
A map between Borel spaces is \definedterm{Borel} if preimages of Borel sets are Borel.

Let $T \from X \to X$ be a Borel map.
We say that a set $Y \subseteq X$ is \definedterm{forward $T$-invariant}\footnote{Also called \definedterm{$T$-absorbing} in the literature.} if $T (Y) \subseteq Y$, and \definedterm{$T$-complete} if it meets every $T$-orbit.

We approximate null sets, or sets of arbitrarily small measure, as follows. 
A sequence $\sequence{Z_i}[i \in \N]$ is \definedterm{vanishing} if it is decreasing and $\intersection[i \in \N][Z_i] = \emptyset$. 
Indeed, if the sets $Z_i$ were Borel, then the integrals $\int_{Z_i} |f| \dm$, for integrable functions $f$, would be arbitrarily small for large enough $i$.

Turning now to the Borel translation of the existence of the limit $\limit{f}{g}$ a.e., we let 
\[
\limitsuperior{f}{g} \defeq \limsup_{n \to \infty} \ratio{f}{g}{n}
\]
and note that for each $x \in X$, the limit $\limit{f}{g}(x)$ exists if for every $\e > 0$, we have $\ratio{f}{g}{n}(x) \ge \limitsuperior{f}{g}(x) - \e$ for all large enough $n \in \N$.
To witness the latter condition in a Borel fashion, we need the following notion.
A set $B \subseteq X$ is called \definedterm{uniformly $T$-syndetic} if there is $N \in \N$ such that for each $x \in X$, the gaps in the set $\set{n \in \N : T^n(x) \in B}$ are bounded above by $N$; in other words, $X = \bigcup_{n < N} T^{-n}(B)$ for some $N \in \N$.
We denote by $\hittingtime{T}{B}(x)$ the \textbf{first hitting time} for $B$, i.e.\ the smallest $n \in \positiveintegers \defeq \N \setminus \set{0}$ with $T^n(x) \in B$.
We omit the subscript $T$ from $\hittingtime{T}{B}$ if it is clear from the context.

Our Borel theorem roughly states that for each $\e > 0$, there is a uniformly $T$-syndetic Borel set $B \subseteq X$ such that the 99\% of the interval $\set{x, T x, \dots, T^{n_B(x)-1} x}$ can be tiled with smaller intervals on which the average-ratio of $f$ and $g$ is $\ge \limitsuperior{f}{g}(x) - \e$, as desired.
We then deduce the measurable version (\cref{intro:Dowker_with_limit}) from this via our local--global bridge lemma (\cref{local-global}), which establishes a connection between the local weighted sums and global mean.

\begin{theorem}[Measureless ratio ergodic] \label{intro:measureless}
Let $X$ be a Borel space, $T \from X \to X$ be a Borel transformation, and $\coc \from \N \times X \to \positivereals$ be a Borel $T$-cocycle.
Let $f, g \from X \to \R$ be Borel functions with $g > 0$.
Let $R \from X \to \R$ be a $T$-invariant function with $R < \limitsuperior{f}{g}$.

Then there are a forward $T$-invariant $T$-complete Borel set $X' \subseteq X$, a vanishing sequence $\sequence{Z_i}[i \in \N]$ and a decreasing sequence $\sequence{B_i}[i \in \N]$ of Borel subsets of $X'$ such that for each $i \in \N$, the set $B_i$ is uniformly $T$-syndetic and for all $x \in B_i \setminus Z_i$,
\[
\ratio{f}{g}{n(x)}(x) > R(x),
\]
where $n(x) \defeq \hittingtime{T}{B_i \cup Z_i}(x)$.
Moreover, if $T$ is aperiodic and separable (\cref{sec:decomposition}) then $X' = X$.
\end{theorem}

\subsubsection*{Organization}

In \cref{sec:decomposition}, we establish an elementary decomposition result allowing us to only consider two orthogonal cases: when $T$ is aperiodic and satisfies a local notion of separability, and when $f \times g \times \coc_1$ is eventually periodic along the forward orbits of $T$.
\cref{sec:average-ratios} establishes basic properties of average-ratios used throughout the paper.
In \cref{sec:aperiodic,sec:periodic}, we prove strengthenings of \cref{intro:measureless} in both of the said orthogonal cases. 
In \cref{sec:aperiodic}, we also provide a vanishing markers lemma (\cref{aperiodic:markers}) for an aperiodic separable Borel transformation, which is used in the proof of \cref{intro:measureless} in the aperiodic case. 
We combine the two orthogonal cases of \cref{intro:measureless} in \cref{sec:general}, proving a technical refinement (\cref{measureless_refined}) of \cref{intro:measureless} in the general case.
Finally, in \cref{sec:bridge-Dowker}, we prove our local--global bridge lemma (\cref{local-global}), using which we derive the strengthening of \Dowker's theorem (\cref{intro:Dowker_with_limit,Dowker_with_limit}) from \cref{measureless_refined}.

\subsubsection*{Acknowledgments}

We thank the anonymous referees for their very helpful suggestions and corrections, which significantly improved the paper.
B.D.~Miller was partially supported by FWF grant P29999.
A.~Tserunyan was partially supported by NSF grant DMS-1855648 and NSERC Discovery Grant RGPIN-2020-07120.

\section{Decomposition} \label{sec:decomposition}

Throughout, let $X$ be a Borel space and $T \from X \to X$ be a Borel map.

Given a binary relation $R$ on $X$, we say that a family $\calB$ of
subsets of $X$ \definedterm{separates $R$-related points} if, for all distinct
$x \mathrel{R} y$, there exists $B \in \calB$ such that $x \in B$ but $y
\notin B$. 
We say that $R$ is \definedterm{separable} if there is a countable family of Borel sets that separates $R$-related points. 
We also say that a Borel space $Y$ is \definedterm{separable} if $R \defeq Y^2$ is separable, i.e.\ there is a countable collection $\set{B_n}_{n \in \N}$ of Borel subsets of $Y$ \textbf{separating points}, that is: for all distinct $y,y' \in X$ there is $n \in \N$ such that $y \in B_n$ but $y' \notin B_n$.

For $p \in \positiveintegers$, the \definedterm{$T$-period $p$ part} of $h \from X \to Y$ is given by 
\[
\periodicpart[p]{T}{h}
\defeq 
\set{x \in X}[\forall k \in \N \; h(T^k x) = h(T^{k+p} x)].
\] 
In particular, $\periodicpart[p]{T}{h}$ contains the set of fixed points of $T^p$.
If $X$ and $Y$ are Borel spaces, $Y$ is separable, and $h$ and $T$ are
Borel, then the fact that the class of Borel functions is closed under
appropriate compositions and products ensures that $\periodicpart[p]{T}{h}$
is Borel.
Call the set 
\[
\periodicpart{T}{h} \defeq \union[p \in \positiveintegers][{\periodicpart[p]{T}{h}}]
\]
the \definedterm{$T$-periodic part} of $h$.

\begin{lemma} \label{decomposition:periodp}
Let $T \from X \to X$ be a Borel map and $p \in \positiveintegers$.
If $\periodicpart[p]{T}{h} = \emptyset$ for some Borel function $h \from X \to Y$ to some separable Borel space $Y$, then the graph of $T^p$ is separable.
\end{lemma}

\begin{lemmaproof}
Fix a countable family $\calB$ of Borel subsets of $Y$ that separates
points. We need only show that, for all $x \in X$, there exist $B \in \calB$
and $k \in \N$ with $x \in \preimage{(h \composition T^k)}{B}$ but $T^p(x)
\notin \preimage{(h \composition T^k)}{B}$. Since $x \notin \periodicpart[p]{T}{h}$, so there exists $k \in \N$ such that $(h \composition T^k)(x) \neq (h
\composition T^{k+p})(x)$, and thus there exists $B \in \calB$ such that $(h \composition T^k)(x) \in B$ but $(h \composition T^{k+p})(x) \notin B$.
\end{lemmaproof}

We say that $T \from X \to X$ is \definedterm{aperiodic} if its positive powers are fixed-point free. 
We say that $T$ is \definedterm{separable} if its positive powers have separable graphs. \cref{decomposition:periodp} implies the following:

\begin{lemma} \label{decomposition:aperiodic}
Let $T \from X \to X$ be a Borel map.
If $\periodicpart{T}{h} = \emptyset$ for some Borel function $h \from X \to Y$ to some separable Borel space $Y$, then $T$ is aperiodic and separable.
\end{lemma}

For a set $Y \subseteq X$, we denote by $\backwardorbit{Y}{T} \defeq \bigcup_{n \in \positiveintegers} T^{-n} Y$.
Note that $\backwardorbit{Y}{T}$ is Borel if $Y$ is, and $\backwardorbit{Y}{T}$ is $T$-invariant if $Y$ is forward $T$-invariant.

We say that a function $h \from X \to Y$ is \definedterm{eventually $T$-periodic} if $\periodicpart{T}{h}$ is $T$-complete.; equivalently, $\backwardorbit{\periodicpart{T}{h}}{T} = X$.
Call a partition $X = \Xas \cup \Xep$ \textbf{$h$-admissible for $T$} if $\Xas$ and $\Xep$ are $T$-invariant Borel sets such that $\restriction{T}{\Xas}$ is aperiodic and separable, and $\restriction{h}{\Xep}$ is eventually $\restriction{T}{\Xep}$-periodic.
With this terminology, \cref{decomposition:aperiodic} gives the following:

\begin{proposition}[Decomposition]\label{decomposition}
Let $X$ be a Borel space and let $T \from X \to X$ be a Borel map.
Let $h \from X \to Y$ be a Borel function to some separable Borel space $Y$.
Then there is an $h$-admissible partition $X = \Xas \cup \Xep$ for $T$, namely $\Xep \defeq \backwardorbit{\periodicpart{T}{h}}{T}$.
\end{proposition}

This decomposition allows us to handle the two corresponding cases of \cref{intro:measureless} separately, where we take $h \defeq f \times g \times \coc_1 \from X^3 \to \R \times \positivereals \times \positivereals$ in the notation of the theorem.

\section{Cocycle-weighted average-ratios}\label{sec:average-ratios}

Throughout this section, we let $X$ be a set, $T \from X \to X$ be a map, and $\coc \from \N \times X \to \positivereals$ be a $T$-cocycle. 
We also let $f,g \from X \to \R$ be functions such that $g > 0$.

\begin{proposition}[Convexity of ratios]\label{convexity_of_ratios}
For all $m,n \in \N$ and $x \in X$, we have 
\[
\ratio{f}{g}{m + n}(x) =  \alpha \cdot \ratio{f}{g}{m}(x) + (1-\alpha) \cdot \ratio{f}{g}{n}(T^m x),
\]
where $\alpha \in [0,1]$, namely, $\alpha \defeq \frac{\uppersum{g}{m}}{\uppersum{g}{m + n}}$.
\end{proposition}
\begin{propositionproof}
By the cocycle identity \labelcref{eq:coc_identity}, $\coc_{m+k}(x) = \coc_m(x) \cdot \coc_k(T^m x)$ for all $k \in \N$, so for every function $h \from X \to \R$ and $x \in X$, we have $\uppersum{h}{m+n}(x) = \uppersum{h}{m}(x) + \coc_m(x) \cdot \uppersum{h}{n}(T^m x)$.
Thus,
\begin{align*}
\ratio{f}{g}{m+n}(x)
&= 
\frac{\uppersum{f}{m}(x)}{\uppersum{g}{m+n}(x)} + \frac{\coc_m(x) \cdot \uppersum{f}{n}(T^m x)}{\uppersum{g}{m+n}(x)}
\\
&= 
\frac{\uppersum{g}{m}(x)}{\uppersum{g}{m + n}(x)} \cdot \ratio{f}{g}{m}(x) + \frac{\coc_m(x) \cdot \uppersum{g}{n}(T^m x)}{\uppersum{g}{m + n}(x)} \cdot \ratio{f}{g}{n}(T^m x)
\\
&= 
\alpha \cdot \ratio{f}{g}{m}(x) + (1-\alpha) \cdot \ratio{f}{g}{n}(T^m x).
\qedhere
\end{align*}
\end{propositionproof}

For $x \in X$, let $\limitpoints{f}{g}(x)$ denote the set of limit points of $\set{\ratio{f}{g}{n}(x)}[n \in \N]$.

\begin{proposition} \label{lim_invariance}
If $\lim_n \uppersum{g}{n}(x) = \infty$ for all $x \in X$, then $\limitpoints{f}{g}$ is $T$-invariant.
In particular, the functions $\limitsuperior{f}{g}$ and $\limitinferior{f}{g}$ are $T$-invariant.
\end{proposition}

\begin{propositionproof}
For each $x \in X$ and $n \in \N$, \cref{convexity_of_ratios} gives
\[
\frac{g(x)}{\uppersum{g}{n+1}(x)} \cdot \frac{f(x)}{g(x)} + \frac{\uppersum{g}{n+1}(x) - g(x)}{\uppersum{g}{n+1}(x)} \cdot \ratio{f}{g}{n}(T x).
\]
The first summand converges to $0$, and $\frac{\uppersum{g}{n+1}(x) - g(x)}{\uppersum{g}{n+1}(x)}$ converges to $1$ as $n \to \infty$, so along any subsequence $(n_k)$, we have $\lim_k \ratio{f}{g}{n_k + 1}(x) = \lim_k \ratio{f}{g}{n_k}(T x)$, which shows that $\limitpoints{f}{g}(x) = \limitpoints{f}{g}(T x)$.
\end{propositionproof}

\begin{lemma}\label{periodic:coc_invariant}
    $\coc_p$ is $T$-invariant on $\periodicpart[p]{T}{\coc_1}$ for each $p \in \positiveintegers$, i.e.\ $\coc_p(x) = \coc_p(T x)$ for every $x \in \periodicpart[p]{T}{\coc_1}$.
\end{lemma}
\begin{lemmaproof}
    By the cocycle identity \labelcref{eq:coc_identity}, 
    \[
    \coc_p(x) = \prod_{i < p} \coc_1(T^i x) = \coc_1(x) \cdot \prod_{1 \le i < p} \coc_1(T^i x) = \prod_{1 \le 1 \le p} \coc_1(T^i x) = \coc_p(T x).
    \qedhere
    \]
\end{lemmaproof}

\begin{lemma}\label{coc_ge_1}
    For each $p \in \positiveintegers$ and $x \in \periodicpart[p]{T}{g \times \coc_1}$, if $\lim_n \uppersum{g}{n}(x) = \infty$ then $\coc_p(x) \ge 1$.
\end{lemma}
\begin{lemmaproof}
    The cocycle identity \labelcref{eq:coc_identity} implies that $x \in \periodicpart[p]{T}{\coc_p}$, and hence $\coc_{ip}(x) = \coc_p(x)^i$ for all $i \in \N$.
    Furthermore, the cocycle identity again and period $p$ give
    \[
    \uppersum{g}{n p}(x) = \sum_{i < n} \coc_{ip}(x) \cdot \uppersum{g}{p}(T^{ip} x) = \sum_{i < n} \coc_p(x)^i \cdot \uppersum{g}{p}(x),
    \]
    which is bounded above by $\frac{\uppersum{g}{p}(x)}{1 - \coc_p(x)}$ if $\coc_p(x) < 1$, contradicting $\lim_n \uppersum{g}{n p}(x) = \infty$.
\end{lemmaproof}

\section{The aperiodic separable case} \label{sec:aperiodic}

For $T \from X \to X$, $Y \subseteq X$, and $i \in \N$, put \heightcorrection{$\preimage[\le i]{T}{Y} \defeq \union[j \le i][{\preimage[j]{T}
{Y}}]$}.
We denote by $\setcomplement Y$ the complement of $Y$, i.e.\ $X \setminus Y$.

The proof of the marker lemma for Borel automorphisms (see \cite[Lemma 1
of \S3]{SlamanSteel}) generalizes beyond the injective and standard cases, yielding the following characterization of aperiodicity for separable Borel transformations.

\begin{proposition}[Uniformly syndetic markers] \label{aperiodic:markers}
Let $X$ be a Borel space and $T \from X \to X$ be a separable Borel transformation. 
The following are equivalent:
\begin{enumerate}[(1)]
    \item \label{part:aperiodic} $T$ is aperiodic.

    \item \label{part:markers} There is a vanishing sequence $\sequence{B_i}[i \in \N]$ of $T$-complete Borel subsets of $X$.

    \item \label{part:bdd-markers} There is a vanishing sequence
$\sequence{B_i}[i \in \N]$ of uniformly $T$-syndetic Borel subsets of $X$.
\end{enumerate}
\end{proposition}

\begin{propositionproof} The implications \labelcref{part:bdd-markers}$\imp$\labelcref{part:markers}$\imp$\labelcref{part:aperiodic} are immediate/trivial.

For \labelcref{part:aperiodic}$\imp$\labelcref{part:markers}, assume that $T$ is aperiodic and set $R \defeq \union[n \in \positiveintegers][\graph{T^n}]$. Then the separability
of $T$ yields a family $\set{A_i}[i \in \N]$ of Borel subsets of $X$ that
separates $R$-related points. Put \heightcorrection{$A_s \defeq \intersection[i
\in \preimage{s}{\set{0}}][A_i] \intersection \intersection[i \in \preimage{s}
{\set{1}}][\protect \setcomplement{A_i}]$} for all $s \in \Cantortree$, let
$\le$ denote the lexicographical ordering of $\Cantortree$, and define 
\[
s_i(x) \defeq \textstyle \min_{\le} \set{s \in \Cantorspace[i]}[A_s \intersection \forwardorbit{x}{T} \text{ is infinite}]
\] 
for all $i \in \N$ and $x \in X$, where $\forwardorbit{x}{T} \defeq \set{x, Tx, T^2 x, \dots}$ is the forward $T$-orbit of $x$. As $s_i$ is $T$-invariant, the intersection of the set 
\[
B_i' \defeq \set{x \in X}[x \in A_{s_i(x)}] = \bigcup_{s \in \Cantortree[i]} A_s \intersection \preimage{s_i}{\set{s}}
\]
with each forward orbit of $T$ is infinite. 
The fact that $A_s = A_{s
\concatenation \sequence{0}} \disjointunion A_{s \concatenation \sequence
{1}}$ for all $s \in \Cantortree$ ensures that \heightcorrection{$\preimage
{s_i}{\set{s}} = \preimage{s_{i+1}}{\set{s \concatenation \sequence{0}}}
\disjointunion \preimage{s_{i+1}}{\set{s \concatenation \sequence{1}}}$} for
all $i \in \N$ and $s \in \Cantortree[i]$, so $B_{i+1}' \subseteq B_i'$ for all $i
\in \N$. The set $B' \defeq \intersection[i \in \N][B_i']$ intersects each forward
orbit of $T$ in at most one point, for if $k \in \positiveintegers$ and $x \in
B'$, then there exists $i \in \N$ such that $x \in A_i$ and $T^k(x) \notin
A_i$, so the fact that $x \in B_{i+1}'$ implies that $T^k(x) \notin B_{i+1}'$,
and thus $T^k(x) \notin B'$. The sets $B_i \defeq B_i' \setminus B'$ are therefore as desired.

For \labelcref{part:markers}$\imp$\labelcref{part:bdd-markers}, let $\sequence{A_i}[i \in \N]$ be a vanishing sequence of $T$-complete Borel sets with $A_0 = X$. 
For all $i \in \N$, define 
$
B_i \defeq
A_i \union \union[j < i][{A_j \setminus \preimage[\le i]{T}{A_{j+1}}}]
$.
To see that $B_{i+1} \subseteq B_i$ for all $i \in \N$, note that $A_{i+1},
A_i \setminus \preimage[\le i+1]{T}{A_{i+1}} \subseteq A_i$ and $A_j
\setminus \preimage[\le i+1]{T}{A_{j+1}} \subseteq A_j \setminus
\preimage[\le i]{T}{A_{j+1}}$ for all $j < i$. To see that $\intersection[i
\in \N][B_i] = \emptyset$, note that if $j \in \N$ and $x \in A_j \setminus
A_{j+1}$, then there exists $i > j$ for which $x \in \preimage[\le i]{T}
{A_{j+1}}$, so $x \notin B_i$. Finally, to see that $X = \preimage[\le i^2]{T}
{B_i}$ for all $i \in \N$, recall that $X = A_0$ and note that if $j < i$, then $A_j \subseteq B_i \union
\preimage[\le i]{T}{A_{j+1}}$, so $A_j \subseteq \preimage[\le i(i-j)]{T}{B_i}$
by induction on $i - j$. Thus the sequence $\sequence{B_i}[i \in \N]$ is as desired.
\end{propositionproof}

We now prove a refinement of \cref{intro:measureless} for aperiodic and separable transformations.

\begin{proposition}\label{aperiodic:main}
Let $X$ be a Borel space, $T \from X \to X$ be an aperiodic and separable Borel transformation, and $\coc \from \N \times X \to \positivereals$ be a Borel $T$-cocycle.
Let $f, g \from X \to \R$ be Borel functions with $g > 0$.
Let $R \from X \to \R$ be a $T$-invariant function with $R < \limitsuperior{f}{g}$.

Then there are vanishing sequences $\sequence{Z_i}[i \in \N]$ and $\sequence{B_i}[i \in \N]$ of Borel subsets of $X$ such that for each $i \in \N$, the set $B_i$ is uniformly $T$-syndetic and for all $x \in X \setminus Z_i$,
\[
\ratio{f}{g}{n(x)}(x) > R(x),
\]
where $n(x) \defeq \hittingtime{T}{B_i \cup Z_i}(x)$.
\end{proposition}
\begin{propositionproof}
By \cref{aperiodic:markers}, there is a vanishing sequence $\sequence{B_i}[i \in \N]$ of uniformly $T$-syndetic Borel sets.
Then the sequence given by 
\[
Z_i \defeq \set{x \in X}[\ratio{f}{g}{n(x)}(x) \le R(x) \text{ for all } n < \hittingtime{T}{B_i}(x)]
\]
is also vanishing.
Now fix $i \in \N$ and for each $x \in X \setminus Z_i$, let $m(x) \in \N$ be the least number such that $\ratio{f}{g}{m(x)}(x) > R(x)$.
Note that for each $x \in X \setminus Z_i$, we have $m(x) \le \hittingtime{T}{B_i}(x)$ by the definition of $Z_i$. 
In fact, we claim that $m(x) \le n(x) \defeq \hittingtime{T}{B_i \cup Z_i}(x)$.
Indeed, if $T^k x \in Z_i$ for some $k < m(x)$, then $\hittingtime{T}{B_i}(T^k x) = m(x)-k$, so $\ratio{f}{g}{m(x) - k}(T^k x) \le R(T^k x) = R(x)$, which implies $\ratio{f}{g}{k}(x) > R(x)$ because of the convexity of ratios (\cref{convexity_of_ratios}), namely,
\[
\ratio{f}{g}{m(x)}(x) = \alpha \cdot \ratio{f}{g}{k}(x) + (1 - \alpha) \cdot \ratio{f}{g}{m(x) - k}(T^k x),
\]
for some $\alpha \in (0,1)$, contradicting the minimality of $m(x)$.

Now fix $x \in X \setminus Z_i$. 
Inductively define a sequence $0 = k_0 < k_1 < \dots < k_\ell = n(x)$ of natural numbers by setting $k_0 \defeq 0$ and $k_{i+1} \defeq m(T^{k_i} x)$ if $k_i < n(x)$.
Then observe that the set $\set{x, T x, \dots, T^{n(x) - 1} x}$ is partitioned into sets (tiles) of the form $\set{T^{k_i} x, T^{k_i + 1} x, \dots, T^{k_{i+1} - 1} x}$ and the ratio over each tile is $> R(x)$; more precisely, 
\[
\ratio{f}{g}{k_{i+1} - k_i}(T^{k_i} x) > R(T^{k_i} x) = R(x).
\]
It now follows from the convexity of ratios (\cref{convexity_of_ratios}) that $\ratio{f}{g}{n(x)} > R(x)$.
\end{propositionproof}

\section{The eventually periodic case} \label{sec:periodic}

The set of limit points of $\set{\ratio{f}{g}{n}(x)}[n \in \N]$ is easily computed in the periodic case:

\begin{lemma} \label{periodic:exact}
Let $X$ be a set, $T \from X \to X$ be a map, and $\coc \from \N \times X \to \positivereals$ be a $T$-cocycle.
Let $f, g \from X \to \R$ be functions where $g > 0$, and fix $p, r \in \N$ with $p > 0$ and $r < p$.

Then for each \heightcorrection{$x \in \periodicpart[p]{T}{f \times g \times \coc_1}$}, we have:
\begin{enumerate}[(i),leftmargin=20pt]
\item If $\coc_p(x) \le 1$, then $\lim_n \ratio{f}{g}{p n + r}(x) = \ratio{f}{g}{p}(x)$.

\item If $\coc_p(x) > 1$, then $\lim_n \ratio{f}{g}{p n + r}(x) =  \ratio{f}{g}{p} (T^r x)$.
\end{enumerate}
\end{lemma}

\begin{lemmaproof}
For each $k \in \positiveintegers$, the cocycle identity \labelcref{eq:coc_identity} yields $\coc_k = \prod_{i=0}^{k-1} (\coc_1 \circ T^i)$, so $\periodicpart[p]{T}{f \times g \times \coc_1} \subseteq \periodicpart[p]{T}{f \times g \times \coc_k}$.
Therefore, a further application of the cocycle identity \labelcref{eq:coc_identity} ensures that for each $q \in \N$,
\[
\coc_{p n + q}(x) = \coc_q(T^{p n} x) \product[i < n][\coc_p(T^{p i} x)] =
\coc_q(x) \coc_p(x)^n.
\]
It follows that for any function $h \from X \to \R$,
\begin{align*}
%\longalign
\uppersum{h}{p n + r}(x) = 
& 
\sum_{i < n} \summation[q < p][h (T^{p i + q} x) \coc_{p i + q}(x)] \mathrel{+} 
\\
&+ 
\summation[q < r][h (T^{p n + q} x) \coc_{p n + q}(x)] 
\\
=& 
\summation[i < n][{\coc_p(x)^i \summation[q < p][h (T^q x) \coc_q(x)]}]
\\
&+ 
\coc_p(x)^n \cdot \summation[q < r][h (T^q x) \coc_q(x)] 
\\
=& 
W_n \cdot \uppersum{h}{p}(x) + \coc_p(x)^n \cdot \uppersum{h}{r}(x),
\end{align*}
where $W_n \defeq \summation[i < n][\coc_p(x)^i]$.

\smallskip

\noindent \emph{Case 1:} $\coc_p(x) < 1$. 
Then $\lim_n W_n = 1 / (1 - \coc_p(x))$ and $\lim_n \coc_p(x)^n = 0$, so 
\[
\lim_n \uppersum{h}{p n + r}(x) = \uppersum{h}{p}(x) / (1 - \coc_p(x)), 
\]
and hence $\lim_n \ratio{f}{g}{p n + r}(x) = \ratio{f}{g}{p}(x)$.

\smallskip

\noindent \emph{Case 2:} $\coc_p(x) = 1$.
Then $W_n = n$ and $\coc_p(x)^n = 1$, so $\uppersum{h}{p n + r}(x) = n \uppersum{h}{p}(x) + \uppersum{h}{r}(x)$, in which case $\lim_n \ratio{f}{g}{p n + r}(x) = \ratio{f}{g}{p}(x)$.

\smallskip

\noindent \emph{Case 3:} $\coc_p(x) > 1$.
Then $W_n = \frac{\coc_p(x)^n - 1}{\coc_p(x) - 1}$, so $\coc_p(x)^n = W_n (\coc_p(x) - 1) + 1$, and hence
\begin{align*}
\uppersum{h}{p n + r}(x)
&= 
W_n \cdot \big(\uppersum{h}{p}(x) + (\coc_p(x) - 1) \uppersum{h}{r}(x)\big) + \uppersum{h}{r}(x) 
\\
&= 
W_n \cdot \uppersum{h}{p} (T^r x) + \uppersum{h}{r}(x).
\end{align*}
Since $\lim_n W_n = \infty$, we have $\lim_n \ratio{f}{g}{p n + r}(x) = \ratio{f}{g}{p}(T^r x)$.
\end{lemmaproof}

We now derive a refinement of \cref{intro:measureless} when $f \times g \times \coc_1$ is eventually $T$-periodic.

\begin{proposition} \label{periodic:main}
Let $X$ be a Borel space, $T \from X \to X$ be a Borel map, $\coc \from \N \times X \to \positivereals$ be a Borel $T$-cocycle, and $f, g \from X \to \R$ be Borel functions with $g > 0$. 
Suppose that $f \times g \times \coc_1$ is eventually $T$-periodic and put $X' \defeq \periodicpart{T}{f \times g \times \coc_1}$. 

Then there is a uniformly $(\restriction{T}{X'})$-syndetic Borel set $B \subseteq X'$ such that for all $x \in B$,
\[
\ratio{f}{g}{n(x)}(x) = \limitsuperior{f}{g}(x),
\]
where $n(x) \defeq \hittingtime{T}{B}(x)$.
\end{proposition}

\begin{propositionproof}
Define \heightcorrection{$X' \defeq \union[n \in \positiveintegers][{\periodicpart[n]{T}{f \times g \times \coc_1}}]$} and $n \from X' \to \positiveintegers$ by
\[
n(x) \defeq \min \set{n \in \positiveintegers}[{x \in \periodicpart[n]{T}{f \times g \times \coc_1}}]. 
\]
By \cref{periodic:exact}, the set
\[
A \defeq \set{x \in X'}[\ratio{f}{g}{n(x)}(x) = \limitsuperior
{f}{g}(x)]
\]
is $T$-complete.
Endowing the set $Y \defeq \R \times \positivereals \times \positivereals$ with the lexicographical ordering, let $\lexeq$ denote the corresponding lexicographical ordering of $\functions {\N}{Y}$ and define $\phi \from X' \to \functions{\N}{Y}$ by $\phi(x)(k) \defeq (f \times g \times \coc_1)(T^k x)$. 
Then, putting $\forwardorbit{x}{T} \defeq \set{x, Tx, T^2 x, \dots}$, the set 
\[
B \defeq \set{x \in A}[\phi(x) \lexeq \phi(y) \text{ for all } y \in A \intersection \forwardorbit{x}{T}]
\]
is as desired, since $n(x) = \hittingtime{T}{B}(x)$ for all $x \in B$.
\end{propositionproof}

\section{The general case}\label{sec:general}

Here we state and prove a refinement of \cref{intro:measureless} in the general case.

\begin{theorem}\label{measureless_refined}
Let $X$ be a Borel space, $T \from X \to X$ be a Borel transformation, and $\coc \from \N \times X \to \positivereals$ be a Borel $T$-cocycle.
Let $f, g \from X \to \R$ be Borel functions with $g > 0$.
Let $R \from X \to \R$ be a $T$-invariant function with $R < \limitsuperior{f}{g}$.
Put $h \defeq f \times g \times \coc_1$

Then there is an $h$-admissible partition $X = \Xas \cup \Xep$ for $T$, and for every such partition, there are a vanishing sequence $\sequence{Z_i}[i \in \N]$ of Borel subsets of $\Xas$ and a decreasing sequence $\sequence{B_i}[i \in \N]$ of uniformly $T$-syndetic Borel subsets of $X' \defeq \Xas \cup \big(\periodicpart{T}{h} \cap \Xep)$ such that for each $i \in \N$ and $x \in B_i \setminus Z_i$,
\[
\ratio{f}{g}{n(x)}(x) > R(x),
\]
where $n(x) \defeq \hittingtime{T}{B_i \cup Z_i}(x)$.
\end{theorem}

\begin{theoremproof}
An $h$-admissible partition $X = \Xas \cup \Xep$ is given by \cref{decomposition}.
For every such partition, let $(Z_i)_{i \in \N}$ and $(B_i)_{i \in \N}$ be sequences of subsets of $\Xas$ satisfying the conclusion of \cref{aperiodic:main} for $\restriction{T}{\Xas}$.
Likewise, let $B$ be a subset of $\periodicpart{T}{h} \cap \Xep$ satisfying the conclusion of \cref{periodic:main} for $\restriction{T}{\Xep}$.
Then the sequences $(Z_i)_{i \in \N}$ and $(B_i \cup B)_{i \in \N}$ satisfy the conclusion of the present theorem.
\end{theoremproof}

\section{The strengthening of Dowker's theorem from the measureless version} \label{sec:bridge-Dowker}

In this section we derive our strengthening of \Dowker's ratio ergodic theorem (\cref{intro:Dowker_with_limit,Dowker_with_limit}) from \cref{intro:measureless}.
The core of this derivation is the following lemma, which establishes a connection between the local weighted sums and global mean.

Throughout, let $T \from X \to X$ be a positively nonsingular Borel transformation of a $\sigma$-finite measure space $(X,\mu)$, and let $\coc \from \N \times X \to \positivereals$ be the associated Radon--Nikodym cocycle.

\begin{lemma}[local--global bridge] \label{local-global}
Let $T \from X \to X$ be a positively nonsingular Borel transformation of a $\sigma$-finite measure space $(X,\mu)$, and let $\coc \from \N \times X \to \positivereals$ be the associated Radon--Nikodym cocycle.
Then for each $f \in L^1(X,\mu)$ and each uniformly $T$-syndetic Borel set $B \subseteq X$,
\[
\integral{f}{\mu} = \integral[B]{\uppersum{f}{\hittingtime{T}{B}(x)}(x)}{\mu}[x].
\]
\end{lemma}

\begin{propositionproof}
For all $n \in \N$, put $B_n \defeq \union[1 \le k \le n][{\preimage[k]{T}{B}}]$.
Then for each $n \in \N$,
\begin{align*}
\longalign
\integral[\setcomplement{(B \union B_n)}]{f(T^n x) \coc^x(n)}{\mu}[x] 
  \\
  & = 
  \integral[\setcomplement{\preimage{T}{B \union B_n}}]{f(T^{n+1} x) \coc^{T x}(n) \coc^x(1)}{\mu}[x] 
    \\
  & = \integral[\setcomplement{B_{n+1}}]{f(T^{n+1} x) \coc^x(n+1)}{\mu}[x] 
  \\
  & =  
    \integral[B \setminus B_{n+1}]{f(T^{n+1} x) \coc^x(n+1)}{\mu}[x]
    \\
    & \mspace{20mu} \mathrel{+} 
    \integral[\setcomplement{(B \union B_{n+1})}]{f(T^{n+1} x) \coc^x(n+1)}{\mu}[x]
\end{align*}
by \cref{eq:rho-invariance} and the cocycle identity \labelcref{eq:coc_identity}.
As \heightcorrection{$\integral{f}{\mu} = \integral[B \setminus B_0]{f}{\mu} + \integral[\setcomplement{(B
\union B_0)}]{f}{\mu}$}, induction on $n$ ensures that
\[
\integral{f}{\mu} = 
\summation[k \le n][{\integral[B \setminus B_k]{f(T^k x)\coc^x(k)}{\mu}[x]}]
\mathrel{+} 
\integral[\setcomplement{(B \union B_n)}]{f(T^n x) \coc^x(n)}{\mu}[x].
\]
As $B$ is uniformly $T$-syndetic, there exists $n \in \N$ for which $X = B \union
B_n$, so
\[
 \integral{f}{\mu}
   = \summation[k \le n][{\integral[B \setminus B_k]{f(T^k x)\coc^x(k)}{\mu}[x]}] 
   = \integral[B]{\summation[k \le n][\characteristicfunction
    {\setcomplement{B_k}}(x) f(T^k x) \coc^x(k)]}{\mu}[x],
\]
but the last integrand is \heightcorrection{$\uppersum{f}
{\hittingtime{T}{B}(x)}(x)$} for all $x \in B$.
\end{propositionproof}

\begin{lemma}\label{coc_le_1}
Let $h \from X \to [0,\infty)$ be a $\mu$-integrable function.
Then $\coc_p(x) \le 1$ for a.e.\ $x \in \periodicpart[p]{T}{h \times \coc_1}$, for each $p \in \positiveintegers$.
\end{lemma}
\begin{lemmaproof}
It suffices show that $Y \defeq \set{x \in \periodicpart[p]{T}{g \times \coc_1}}[\coc_p(x) > 1]$ is null.
By \cref{periodic:coc_invariant}, the set $Y$ is forward $T$-invariant, so $T^{-1} Y \supseteq Y$.
This, \cref{eq:rho-invariance}, and the period of $h$ yield
\[
\int_{Y} h \dm 
= 
\int_{T^{-p} Y} (h \circ T^p) \coc_p \dm 
\ge 
\int_{Y} (h \circ T^p) \coc_p \dm
=
\int_{Y} h \coc_p \dm,
\]
which implies that $Y$ is null because $\coc_p > 1$ on $Y$ and $\int_{Y} h \dm < \infty$.
\end{lemmaproof}

We call a set $Y \subseteq X$ \definedterm{$\mu$-almost $T$-invariant} if the symmetric difference of $T^{-1} Y$ and $Y$ is $\mu$-null.
We omit $\mu$ from the terminology if it is clear from the context.

\begin{lemma}\label{almost-periodic}
Let $g \from X \to \positivereals$ be a $\mu$-integrable function and suppose that there is a forward $T$-invariant $T$-complete Borel subset $Y$ of $\periodicpart{T}{g \times \coc_1}$ such that $\coc_p(x) = 1$ for a.e.\ $x \in Y \cap \periodicpart[p]{T}{g \times \coc_1}$, for each $p \in \positiveintegers$.
Then the set $\periodicpart{T}{g \times \coc_1}$ and all of its forward $T$-invariant Borel subsets are almost $T$-invariant.
\end{lemma}

\begin{lemmaproof}
It suffices to fix $p$ and show that every forward $T$-invariant Borel subset $Y_p$ of $Y \cap \periodicpart[p]{T}{g \times \coc_1}$ is almost $T$-invariant.
Using $T^{-p} Y_p \supseteq Y_p$, \cref{eq:rho-invariance}, the period of $g$, and $\restriction{\coc_p}{Y_p} = 1$ a.e., we get
\[
\int_{Y_p} g \dm 
= 
\int_{T^{-p} Y_p} (g \circ T^p) \coc_p \dm 
\ge 
\int_{Y_p} g \coc_p \dm
=
\int_{Y_p} g \dm,
\]
so $(T^{-p} Y_p) \setminus Y_p$ is null because $(g \circ T^p) \coc_p > 0$.
By the forward $T$-invariance of $Y_p$, we have $(T^{-p} Y_p) \supseteq T^{-1} Y_p \supseteq Y_p$, so $T^{-1} Y_p \setminus Y_p$ is also null, and hence $Y_p$ is almost $T$-invariant.
\end{lemmaproof}

%For a Borel transformation $T \from X \to X$ on a $\sigma$-finite measure space $(X,\mu)$, we recall that $T$ is said to be \definedterm{conservative} (with respect to $\mu$) if every $T$-wandering Borel subset of $X$ is null, where a subset $W \subseteq X$ is \definedterm{$T$-wandering} if the sets $T^{-n}(W)$, $n \in \N$, are pairwise disjoint.

\begin{lemma}\label{f_ge_Rg}
Let $f, g \in L^1(X,\mu)$ with $g > 0$ and put $h \defeq f \times g \times \coc_1$.
Let $X = \Xas \cup \Xep$ be an $h$-admissible partition for $T$ such that $\coc_p(x) = 1$ for a.e.\ $x \in \periodicpart[p]{T}{h} \cap \Xep$, for each $p \in \positiveintegers$.
Let $R \from X \to \R$ be a $T$-invariant Borel function such that $R \le \limitsuperior{f}{g}$ a.e.\ and $R g \in L^1(X,\mu)$.
Then $\integral{f}{\mu} \ge \integral{R g}{\mu}$.
\end{lemma}

\begin{lemmaproof}
We need only show that $\integral{f}{\mu} \ge \integral{gh}{\mu} - \e$
for all $\e > 0$, so we can assume that $R < \limitsuperior{f}{g}$. 
By \cref{almost-periodic} applied to $\restriction{T}{\Xep}$ and the forward $T$-invariant $\restriction{T}{\Xep}$-complete subset $Y \defeq \periodicpart{T}{h} \cap \Xep$ of $\periodicpart{T}{g \times \coc_1} \cap \Xep$, the set $Y$ is almost $T$-invariant.
Thus, letting $X', Z_i, B_i \subseteq X$ be as in \cref{measureless_refined} applied to the given partition $X = \Xas \cup \Xep$, it follows that $X'$ is conull.
Since $X'$ is also forward $T$-invariant, we may assume without loss of generality that $X' = X$.

Take $i \in \N$ sufficiently large so that \heightcorrection{$\integral[Z_i]{\absolutevalue{f} + \absolutevalue{R g}}{\mu} \le \e$}, and put $Z \defeq Z_i$ and $B \defeq B_i$. 
Then applications of \cref{local-global} to $f$ and to $|f| \characteristicfunction{Z}$ give:
\begin{align*}
\integral{f}{\mu}
& =
\integral[B]{\uppersum{f}{\hittingtime{T}{B}(x)}(x)}{\mu}[x]
\\
&= 
\integral[B]{\uppersum{f \characteristicfunction{X \setminus Z}}{\hittingtime{T}{B}(x)}(x)}{\mu}[x] 
+ 
\integral[B]{\uppersum{f \characteristicfunction{Z}}{\hittingtime{T}{B}(x)}(x)}{\mu}[x]
\\
&\ge 
\integral[B]{\uppersum{f \characteristicfunction{X \setminus Z}}{\hittingtime{T}{B}(x)}(x)}{\mu}[x] 
- 
\integral[B]{\uppersum{\absolutevalue{f} \characteristicfunction{Z}}{\hittingtime{T}{B}(x)}(x)}{\mu}[x]
\\
&= \integral[B]{\uppersum{f \characteristicfunction{X \setminus Z}}
{\hittingtime{T}{B}(x)}(x)}{\mu}[x] 
- 
\integral[Z]{\absolutevalue{f}}{\mu}.
\end{align*}
Furthermore, putting $K(x) = \set{k < \hittingtime{T}{B}(x)}[T^kx \notin Z \text{ and } (k = 0 \text{ or } T^{k-1} x \in Z)]$ for each $x \in B$, we have:
\begin{align*}
\integral[B]{\uppersum{f \characteristicfunction{X \setminus Z}}{\hittingtime{T}{B}(x)}(x)}{\mu}[x] 
&= 
\int_B \sum_{k \in K(x)} \uppersum{f}{\hittingtime{T}{Z}(T^k x)}(T^k x) \coc_k(x) \ d \mu(x) 
\\
&\ge 
\int_B \sum_{k \in K(x)} R(T^k x) \uppersum{g}{\hittingtime{T}{Z}(T^k x)}(T^k x) \coc_k(x) \ d \mu(x)
\\
&= 
\int_B R(x) \sum_{k \in K(x)} \uppersum{g}{\hittingtime{T}{Z}(T^k x)}(T^k x) \coc_k(x) \ d \mu(x)
\\
&= 
\int_B R(x) \uppersum{g \characteristicfunction{X \setminus Z}}{\hittingtime{T}{B}(x)}(x) \ d \mu(x)
\\
&= 
\int_B \uppersum{R g \characteristicfunction{X \setminus Z}}{\hittingtime{T}{B}(x)}(x) \ d \mu(x)
\\
&\ge 
\integral[B]{\uppersum{R g}{\hittingtime{T}{B}(x)}(x)}{\mu}[x] 
- 
\integral[B]{\uppersum{\absolutevalue{R g} \characteristicfunction{Z}}{\hittingtime{T}{B}(x)}(x)}{\mu}[x] 
\\
& = 
\integral[X]{R g}{\mu} - \integral[Z]{\absolutevalue{R g}}{\mu}
\end{align*}
by the cocycle identity \labelcref{eq:coc_identity}, the $T$-invariance of $R$, and \cref{local-global}. 
Thus,
\begin{align*}
\integral[X]{f}{\mu} \ge \integral[X]{R g}{\mu} - \integral[Z]{\absolutevalue{f}
+ 
\absolutevalue{R g}}{\mu} 
\ge \integral[X]{R g}{\mu} - \e,
\end{align*}
which finishes the proof.
\end{lemmaproof}

Recall that $\Exp{h}$ denotes the conditional expectation of a function $h \in L^1(X,\mu)$ with respect to the $\sigma$-algebra $\calB_T$ of all $T$-invariant Borel subsets of $X$.

\begin{proposition}\label{cond_exp}
    For $f,g \in L^1(X,\mu)$ with $g > 0$, the function $R \defeq \frac{\Exp{f}}{\Exp{g}}$ is the unique (up to null sets) $T$-invariant Borel function satisfying
\begin{equation}\label{eq:cond_exp}
    \int_Y f \dm = \int_Y R g \dm, \hspace{1em} \text{for each } Y \in \calB_Y.
\end{equation}
\end{proposition}
\begin{propositionproof}
    The function $R \defeq \frac{\E(f \mid \calB_T)}{\E(g \mid \calB_T)}$ indeed satisfies \labelcref{eq:cond_exp} because for each $Y \in \calB_T$, 
    \[
    \int_Y f \dm = \int \Exp{f} \dm = \int_Y R \Exp{g} \dm = \int_Y R g \dm.
    \]
    Conversely, if a Borel function $R \from X \to \R$ is $T$-invariant and satisfies \labelcref{eq:cond_exp}, then for each $Y \in \calB_T$,
    \[
    \int_Y f \dm = \int_Y R g \dm = \int_Y R \Exp{g} \dm,
    \]
    so the $T$-invariant function $R \Exp{g}$ fulfils the definition of $\Exp{f}$.
\end{propositionproof}

We are now ready to prove \cref{intro:Dowker_with_limit}, which we restate here as \cref{Dowker_with_limit} for the reader's convenience.

\begin{theorem}\label{Dowker_with_limit}
Let $T$ be a positively nonsingular Borel transformation of a $\sigma$-finite measure space $(X,\mu)$ and let $\coc$ be the associated Radon--Nikodym cocycle. 
Let $f, g \from X \to \R$ be Borel functions such that $f$ is $\mu$-integrable, $g > 0$, and $\lim_n \uppersum{g}{n} = \infty$ a.e.
Then the limit
$
\limit{f}{g} \defeq \lim_n \ratio{f}{g}{n}
$
exists a.e., and is finite and $T$-invariant. 
Moreover:
\begin{enumerate}[(i)]
    \item\label{part:g_integrable} If $g$ is $\mu$-integrable, then $\limit{f}{g} = \frac{\Exp{f}}{\Exp{g}}$ a.e.

    \item\label{part:g_not_integrable} If $\int_Y g \dm = \infty$ for each non-null $Y \in \calB_T$, then $\limit{f}{g} = 0$ a.e.
\end{enumerate}
\end{theorem}

\begin{theoremproof}
Firstly, note that $\limitsuperior{f}{g}$ is $T$-invariant by \cref{lim_invariance}.

Let $p \in \positiveintegers$.
\cref{coc_ge_1} applied to $g$ gives $\coc_p(x) \ge 1$ for a.e.\ $x \in \periodicpart[p]{T}{g \times \coc_1}$.
On the other hand, \cref{coc_le_1} applied to $h \defeq |f|$ gives $\coc_p(x) \le 1$ for a.e.\ $x \in \periodicpart[p]{T}{|f| \times \coc_1}$.
Thus, we have $\coc_p(x) = 1$ for a.e.\ $x \in \periodicpart[p]{T}{f \times g \times \coc_1}$.

\smallskip

\noindent \textit{Proof of \labelcref{part:g_integrable}.}\
Define $\limitinferior{f}{g}(x) \defeq \liminf_{n \goesto \infty} \ratio{f}{g}{n}(x)$ for all $x \in X$. 
To show that $\limit{f}{g}$ exists a.e.\ and is equal to $\frac{\Exp{f}}{\Exp{g}}$, it suffices to show that $\integral[Y]{f}{\mu} \ge \integral[Y]{\limitsuperior{f}{g} \, g}{\mu}$ for each $Y \in \calB_T$, because applying this to $-f$ gives
\begin{align*}
\integral[Y]{f}{\mu}
  & = -\integral[Y]{-f}{\mu} \\
  & \le -\integral[Y]{\limitsuperior{-f}{g} \, g}{\mu}
    \\
  & = -\integral[Y]{-\limitinferior{f}{g} \, g}{\mu} \\
  & = \integral[Y]{\limitinferior{f}{g} \, g}{\mu},
\end{align*}
which implies that the limit $\limit{f}{g}$ exists a.e.\ and satisfies the conditions of \cref{cond_exp}.
Furthermore, by restricting $T$ to $Y$, we may assume without loss of generality that $Y = X$.
Let $X = \Xas \cup \Xep$ be an $(f \times g \times \coc_1)$-admissible partition for $T$ given by \cref{decomposition}.

\begin{claim}
The function $\limitsuperior{f}{g} \, g$ is $\mu$-integrable.
\end{claim}

\begin{lemmaproof}
If $\integral{\absolutevalue{\limitsuperior{f}{g}} \, g}{\mu} = \infty$, then
by the monotone convergence theorem,
\[
\integral{\absolutevalue{f}}{\mu} < \integral
{\min \set{\absolutevalue{\limitsuperior{f}{g}}, r} \, g}{\mu},
\]
for sufficiently large $r \in \positivereals$.
Because $\limitsuperior{f}{g}$ is $T$-invariant, so is $\min \set{\limitsuperior{\absolutevalue{f}}{g},r}$.
Hence an application of \cref{f_ge_Rg} to $\absolutevalue{f}$, $g$, and $\min \set{\limitsuperior{\absolutevalue{f}}{g},r}$ (in places of $f, g, R$) and the partition $X = \Xas \cup \Xep$ yield
\[
\integral{\absolutevalue{f}}{\mu} \ge \integral{\min \set{\limitsuperior{\absolutevalue{f}}{g}, r} \, g}{\mu} \ge \integral{\min \set{\absolutevalue{\limitsuperior{f}{g}}, r} \, g}{\mu},
\] 
a contradiction.
\end{lemmaproof}

As $\limitsuperior{f}{g}$ is $T$-invariant, so is the set $B$ of all points at which $\limitsuperior{f}{g}$ is finite. 
Thus an application of \cref{f_ge_Rg} to $f$, $g$, $R \defeq \limitsuperior{f}{g}\characteristicfunction{B}$, and partition $X = \Xas \cup \Xep$ ensures that \heightcorrection{$\integral{f}{\mu}\ge \integral{\limitsuperior{f}{g} \, g}{\mu}$}.

\smallskip

\noindent \textit{Proof of \labelcref{part:g_not_integrable}.}\
It suffices to show that $\limit{|f|}{g} = 0$ a.e., so by replacing $f$ with $|f|$ we may assume without loss of generality that $f \ge 0$.

Suppose toward a contradiction that the set of points $x \in X$ with $\limitsuperior{f}{g}(x) > 0$ is non-null.
Then for a small enough $r > 0$, the set $Y_r$ of points $x \in X$ with $\limitsuperior{f}{g}(x) > r$ is non-null.
Because $Y_r$ is $T$-invariant and Borel, $\int_Y g \dm = \infty$.
Thus by restricting to $Y_r$ we may assume without loss of generality that $X = Y_r$.

By the $\sigma$-finiteness of $\mu$ there is an increasing sequence $(g_k)_{k \in \N}$ of positive $\mu$-integrable functions converging to $g$, so $\lim_k \int g_k \dm = \int g \dm$ by the monotone convergence theorem.
Furthermore, $\ratio{f}{g_k}{n} \ge \ratio{f}{g}{n}$ for all $n \in \positiveintegers$, so $\limitsuperior{f}{g_k} \ge \limitsuperior{f}{g} > r$.

Fix $k \in \N$ and let $X = \Xas^{(k)} \cup \Xep^{(k)}$ be an $(f \times g_k \times \coc_1)$-admissible partition for $T$ given by \cref{decomposition}.
Then the restriction of $T$ to $\Xas \cup \Xas^{(k)}$ is still aperiodic and separable, so $X = (\Xas \cup \Xas^{(k)}) \cup (\Xep \cap \Xep^{(k)})$ is still an $(f \times g_k \times \coc_1)$-admissible partition for $T$.
Thus, an application of \cref{f_ge_Rg} to this partition and functions $f$, $g_k$, $R \defeq r$ gives
\[
\int f \dm \ge \int r g_k \dm \to \int r g \dm = \infty
\]
as $k \to \infty$, contradicting the $\mu$-integrability of $f$.
\end{theoremproof}

\def\MR#1{}
\bibliographystyle{amsalpha}
\bibliography{bibliography}

\end{document}